\documentclass[11pt]{article}%[10pt]{article}%\documentclass
\usepackage{amsmath,latexsym,epsf}%,mathabx}%
\begin{document}

%% ----------------------------------------------------------------------
\newcommand{\nc}{\newcommand}
\def\PP#1#2#3{{\mathrm{Pres}}^{#1}_{#2}{#3}\setcounter{equation}{0}}
\def\ns{$n$-star}\setcounter{equation}{0}
\def\nt{$n$-tilting}\setcounter{equation}{0}
\def\Ht#1#2#3{{{\mathrm{Hom}}_{#1}({#2},{#3})}\setcounter{equation}{0}}
\def\qp#1{{${(#1)}$-quasi-projective}\setcounter{equation}{0}}
\def\mr#1{{{\mathrm{#1}}}\setcounter{equation}{0}}
\def\mc#1{{{\mathcal{#1}}}\setcounter{equation}{0}}
\def\HD{\mr{Hom}_{\mc{D}}}
\def\AdT{\mr{Add}_{\mc{T}}}
\def\adT{\mr{add}_{\mc{T}}}
\def\Kb{\mc{K}^b(\mr{Proj}R)}
\def\kb{\mc{K}^b(\mc{P}_R)}

%%%%%%%%FROM latexexam.tex
%\theoremstyle{definition}
\newtheorem{Th}{Theorem}[section]%[section]%
\newtheorem{Sc}[Th]{}

\newtheorem{Def}[Th]{Definition}
\newtheorem{Lem}[Th]{Lemma}
\newtheorem{Pro}[Th]{Proposition}
\newtheorem{Cor}[Th]{Corollary}
\newtheorem{Rem}[Th]{Remark}
\newtheorem{Exm}[Th]{Example}

\def\Pf#1{{\noindent\bf Proof}.\setcounter{equation}{0}}
\def\>#1{{ $\Rightarrow$ }\setcounter{equation}{0}}
\def\<>#1{{ $\Leftrightarrow$ }\setcounter{equation}{0}}
\def\bskip#1{{ \vskip 20pt }\setcounter{equation}{0}}
\def\sskip#1{{ \vskip 5pt }\setcounter{equation}{0}}
\def\mskip#1{{ \vskip 10pt }\setcounter{equation}{0}}
\def\bg#1{\begin{#1}\setcounter{equation}{0}}
\def\ed#1{\end{#1}\setcounter{equation}{0}}
\def\KET{T^{^F\bot}\setcounter{equation}{0}}
\def\KEC{C^{\bot}\setcounter{equation}{0}}
%\def\KET{\mr{KerExt}_A^{i\ge 1}(T,-)\setcounter{equation}{0}}
%\def\KEC{\mr{KerExt}_A^{i\ge 1}(C,-)\setcounter{equation}{0}}
%\def\KEC{\mr{KerExt}_A^{i\ge 1}(-,C)\setcounter{equation}{0}}
%%%%%%%%%%
\def\e{\delta\setcounter{equation}{0}}
\def\D{\mathbf{Diff(R)}\setcounter{equation}{0}}
\def\jze{{ \begin{pmatrix} 0 & 0 \\ 1 & 0 \end{pmatrix}}\setcounter{equation}{0}}
\def\hjz#1#2{{ \begin{pmatrix} {#1} & {#2} \end{pmatrix}}\setcounter{equation}{0}}
\def\ljz#1#2{{  \begin{pmatrix} {#1} \\ {#2} \end{pmatrix}}\setcounter{equation}{0}}
\def\jz#1#2#3#4{{  \begin{pmatrix} {#1} & {#2} \\ {#3} & {#4} \end{pmatrix}}\setcounter{equation}{0}}

\def\Er{\mr{Ext}_R}
\def\ED{\mr{Ext}_{\D}}
\def\Hr{\mr{Hom}_R}
\def\HD{\mr{Hom}_{\D}}
\def\RD#1{(#1\oplus #1,\jze)}

%%%%%%%%%%%%%%%%%%%%%%%%%%%%%%%%%%%%%%%%%%%%%%%%%%%%%%%%%%%%%%%%%%%%%%%%%%%%%%%%%
%**************************±êÌâ¡¢ÕªÒª¡¢·ÖÀàºÅ¡¢¹Ø¼ü×Ö**************************

\title{\bf  Gorenstein Homological Theory for Differential Modules {\thanks {Supported by the National Science Foundation of China
(Grant Nos. 10971099 and 11171149)}}}
\smallskip
\author{\small {Jiaqun WEI}\\
\small Institute of Mathematics, School of Mathematics Sciences,\\
\small
Nanjing Normal University %\\
\small Nanjing 210046, P.R.China\\ \small Email:
weijiaqun@njnu.edu.cn}
\date{}
\maketitle
\baselineskip 15pt%16pt%14pt%15.5pt%\baselineskip  25.5pt %
%%%%%%%\hskip 18pt
%
% Abstract ------------------------------------------------------
%
\begin{abstract}
\vskip 10pt%
We show that a differential module is Gorenstein projective if and
only if its underlying module is Gorenstein projective. Dually, a
differential module is Gorenstein injective if and only if its
underlying module is Gorenstein injective.
%
%
%\mskip\
%
%\noindent 2000 Mathematics Subject Classification: Primary 18E30
%16E05 Secondary 18G35 16G10
%
%%\sskip\
%
%\noindent {\it Keywords}: semi-tilting complex, tilting module,
%derived category

%
\end{abstract}
%\smallskip
%
\vskip 30pt
% ----------------------------------------------------------------------
%% ----------------------------------------------------------------------
%\def\baselinestretch{1}

\section{Introduction}
%{\noindent \Large\bf Introduction}

%%
%---------------
%PUSHOUT and PULLBACK
%
%---------------
%

Let $R$ be a ring. A differential module $(M,\e)$ over $R$ is by
definition an $R$-module $M$ equipped with  an $R$-endomorphism
$\e$ of square zero. A homomorphism from a differential module
$(M,\e_M)$ to  a differential module $(N,\e_N)$ is an
$R$-homomorphism from $M$ to $N$ such that $\e_Mf=f\e_N$. We
denote by $\D$ the category of all differential $R$-modules. Note
that it is just the module category of the ring of dual numbers
over $R$ [\ref{CE}]. In particular, it admits enough projective
objects and enough injective objects.

The notion of differential modules already appeared in Cartan and
Eilenberg's book [\ref{CE}] five decades ago. However, it is
recent thing that the study of differential modules is founded
interesting in their own right. In the paper [\ref{ABI}], Avaromov
etc. studied class and rank of differential modules, given common
generalizations of important results in commutative algebra and
graded polynomial rings. Ringel and Zhang [\ref{RZ}] recently
provided interesting relationship between the module category of a
hereditary (artin) algebra $R$ and the stable category of
Frobenius category $\mathcal{L}$ of perfect differential modules
over $R$, and they also described the representation theory of
$\mathcal{L}$.

This paper focuses on Gorenstein homological theory of
differential modules. Gorenstein homological theory originated in
the works of Auslander and Bridger [\ref{AuB}], where they
introduced G-dimensions. Enochs extended their ideas and
introduced Gorenstein projective,  Gorenstein injective and
Gorenstein flat modules and correspondent dimensions over
arbitrary rings, see the book [\ref{Eb}] for details. Later
Gorenstein homogical theory was extensively studied and developed
by Avromov, Martsinkovsky, Christensen, Veliche, Sather-Wagstaff,
Chen, Beligiannis, Yang and many others (see for instance
[\ref{AM}, \ref{Be},
 \ref{Chen}, \ref{Ch}, \ref{Gao}, \ref{SSW}, \ref{V}, \ref{Y}] etc.).

Our main result states as follows.

\bg{Th}\label{MT}%
Let $R$ be a ring and $(M,\e)$ a differential module. Then
$(M,\e)$ is Gorenstein projective (in the category of differential
modules) if and only if $M$ is Gorenstein projective (in the
category of $R$-modules). Dually, $(M,\e)$ is Gorenstein injective
if and only if $M$ is Gorenstein injective.%
\ed{Th}

The theorem provides interesting relationships between Gorestein
homology theory of differential modules and that of their
underlying modules. Furthermore, we have the following result.
Here, $\mr{Gpd}M$ (respectively, $\mr{Gid}M$) means the Gorenstein
projective (respectively, Gorenstein injective) dimension of the
module $M$ and $\mr{Ggd}R$ means the Gorenstein global dimension
of $R$.

\bg{Th}\label{MGdT}%Gorenstein dimension Theorem

Let $(M,\e_M)$ be a differential module and $n$ be an integer.

$(1)$ $\mr{Gpd}(M,\e_M)\le n$  if and only if $\mr{Gpd}M\le n$.

$(2)$ $\mr{Gid}(M,\e_M)\le n$  if and only if $\mr{Gid}M\le n$.

$(3)$ $\mr{Ggd}\D\le n$ if and only if $\mr{Ggd}R\le n$.

\ed{Th}

%%%%%%%%%%%%%%%%%%%%%%%%%%%%%%%%%%%%%%%%%%%%%%%%%%%%%%%%%%%%%%%%%%%%%%%%%
\vskip 30pt
\section{Basic results on differential modules}

We write the element in a direct sum as a row vector. For two maps
$f:M\to N$ and $g: N\to L$, we use $fg$ to denote the composition
of these two maps. Throughout the paper, we fix $R$ a ring (which
is associative with an identity).

%\bg{Sc}\label{TEST}\ed{Sc}\vskip -10pt%projective diff mod The projective

Let $X$ be an $R$-module. It is easy to see that $\RD{X}$ is a
differential module and such differential module is called
contractible [\ref{ABI}].

The following easy lemma is useful in describing homomorphisms
between a differential module and a contractible differential
module. We leave its proof to the reader.

\bg{Lem} \label{mapL} %map Lemma

Let $(M,\e)$ be a differential module and $X$ an $R$-module.

$(1)$ Let $f\in\mr{Hom}_R(M,X\oplus X)$. Then
$f\in\mr{Hom}_{\D}((M,\e),\RD{X})$ if and only if
$f=\begin{pmatrix} g & \e g\end{pmatrix}$ for some
$g\in\mr{Hom}_R(M,X)$.

$(2)$ Let $f\in\mr{Hom}_R(X\oplus X,M)$. Then
$f\in\mr{Hom}_{\D}(\RD{X},(M,\e))$ if and only if
$f=\bg{pmatrix}g\e\\g\ed{pmatrix}$ for some $g\in\mr{Hom}_R(X,M)$.

\ed{Lem}

By the above lemma, we can represent every homomorphism from the
differential module $(M,\e)$ to a contractible differential module
$\RD{X}$ as the form $\bg{pmatrix}g & \e g \ed{pmatrix}$, where
$g\in\mr{Hom}_R(M,X)$. Similarly, we represent every homomorphism
from $\RD{X}$ to $(M,\e)$ as the form $\bg{pmatrix}g\e \\ g
\ed{pmatrix}$.

\bg{Lem} \label{facL} %factor Lemma ·Ö½âÒýÀí

Let $f: (M,\e_M)\to(N,\e_N)$ be a homomorphism between
differential modules and $X$ an $R$-module.

$(1)$ $\mr{Hom}_R(f,X)$ is an epimorphism if and only if
$\mr{Hom}_{\D}(f,\RD{X})$ is an epimorphism.

$(2)$  $\mr{Hom}_R(X,f)$ is an epimorphism if and only if
$\mr{Hom}_{\D}(\RD{X},f)$ is an epimorphism.

\ed{Lem}

\Pf. (1)  Assume first $\mr{Hom}_R(f,X)$ is an epimorphism. Given
a homomorphism $\hjz{g}{\e_Mg}\in \mr{Hom}_{\D}((M,\e_M),\RD{X})$,
where $g\in\mr{Hom}_R(M,X)$, then we have that $g=fh$ for some
$h\in\Hr(N,X)$, by assumptions. Since $\e_Mf=f\e_N$, we see that
$\hjz{g}{\e_Mg}=f\hjz{h}{\e_Nh}$ with $\hjz{h}{\e_Nh}\in
\mr{Hom}_{\D}((N,\e_N),\RD{X})$. It follows that%This shows
$\mr{Hom}_{\D}(f,\RD{X})$ is an epimorphism.

Assume now $\mr{Hom}_{\D}(f,(X\oplus X,\jze))$ is an epimorphism.
Given a homomorphism $g\in\Hr(M,X)$, we have
$\hjz{g}{\e_Mg}\in\HD((M,\e_M),\RD{X})$. Hence there is some
$\hjz{h}{\e_Nh}\in\HD((N,\e_N),\RD{X})$ such that
$\hjz{g}{\e_Mg}=f\hjz{h}{\e_Nh}$, where $h\in\Hr(N,X)$. It follows
that $g=fh$ and so $\mr{Hom}_R(f,X)$ is an epimorphism.

(2) The proof is dual to that of part (1). \hfill$\Box$

%%%%%--------------------------------------------------------------------------

\vskip 10pt

We have the following useful corollary.

\bg{Cor}\label{exC}%exact corollary

Assume that $ 0\to (M,\e_M)\to (N,\e_N)\to (L,\e_L)\to 0\ \
(\dag)$ is an exact sequence of differential modules and $X$ is an
$R$-module. Then

$(1)$ $\HD(\dag, \RD{X})$ is exact if and only if $\Hr(\dag, X)$
is exact.

$(2)$ $\HD(\RD{X},\dag)$ is exact if and only if $\Hr(X,\dag)$ is
exact.

\ed{Cor}

%%%%%--------------------------------------------------------------------------
\vskip 30pt

\section{Gorenstein differential modules}

%%%%%--------------------------------------------------------------------------

We begin with the following general result.

\bg{Lem}\label{pcL}%precover lemma

Let $(M,\e)$ be a differential module and $\mc{C}$ be a class of
$R$-modules. Assume that

\vskip 5pt\hskip 90pt$0\to L\stackrel{\lambda}{\to}
C\stackrel{\pi}{\to} M\to 0\ $\hfill$(\ddag)$\vskip 5pt

\noindent is an exact sequence of $R$-modules such that
$C\in\mc{C}$ and $\Hr(C',\ddag)$ is exact for any $C'\in\mc{C}$.
Then there is an exact sequence of differential modules

\hskip 70pt$0\to (C\oplus L,\e_{C\oplus
L})\stackrel{\jz{-1}{h}{0}{\lambda}}{\longrightarrow}
\RD{C}\stackrel{\ljz{\pi\e}{\pi}}{\to} (M,\e)\to 0,$\hfill
$(\ddag\ddag)$

\noindent for some $h\in\Hr(C,C)$, such that

$(1)$ $\HD(\RD{C'},\ddag\ddag)$ is exact for any $C'\in\mc{C}$.

$(2)$ For any $R$-module $X$, $\Hr(\ddag,X)$ is exact if and only
if $\HD(\ddag\ddag,\RD{X})$ is exact.

\ed{Lem}

\Pf. Since $\Hr(C,\ddag)$ is exact by assumptions, $\Hr(C,\pi)$ is
epi. Hence, there is some $h\in\Hr(C,C)$ such that $\pi\e=h\pi$.
Now we can construct the following commutative diagram.

\bskip\

 \setlength{\unitlength}{0.09in}
 \begin{picture}(50,7)

 \put(12,0){\makebox(0,0)[c]{$0$}}
                             \put(14,0){\vector(1,0){2}}
 \put(18,0){\makebox(0,0)[c]{$L$}}
                             \put(22,1){\makebox(0,0)[c]{$\lambda$}}
                             \put(21,0){\vector(1,0){2}}
 \put(27,0){\makebox(0,0)[c]{$C$}}
                             \put(31,1){\makebox(0,0)[c]{$\pi$}}
                             \put(30,0){\vector(1,0){2}}
 \put(35,0){\makebox(0,0)[c]{$M$}}
                            \put(37,0){\vector(1,0){2}}
 \put(41,0){\makebox(0,0)[c]{$0$}}

                 \put(18.5,4){\line(0,-1){2}}
                 \put(18,4){\line(0,-1){2}}
                 \put(27,4){\vector(0,-1){2}}
                       \put(25,3){\makebox(0,0)[c]{${ \ljz{h}{\lambda}}$}}
                 \put(35,4){\vector(0,-1){2}}
                       \put(37,3){\makebox(0,0)[c]{$\pi\e$}}

 \put(12,6){\makebox(0,0)[c]{$0$}}
                             \put(14,6){\vector(1,0){2}}
 \put(18,6){\makebox(0,0)[c]{$L$}}
                             \put(22,7){\makebox(0,0)[c]{$_{(0\ 1)}$}}
                             \put(21,6){\vector(1,0){2}}
 \put(27,6){\makebox(0,0)[c]{$C\oplus L$}}
                             \put(31,8){\makebox(0,0)[c]{${\ljz{1}{0}}$}}
                             \put(30,6){\vector(1,0){2}}
 \put(35,6){\makebox(0,0)[c]{$C$}}
                             \put(37,6){\vector(1,0){2}}
                             \put(33,5){\vector(-1,-1){4}}
                             \put(32,3){\makebox(0,0)[c]{$h$}}
 \put(41,6){\makebox(0,0)[c]{$0$}}

\end{picture}
\bskip\

It is easy to see that the above diagram is a pullback. Hence
there is an exact sequence

\hskip 70pt$0\to C\oplus
L\stackrel{\jz{-1}{h}{0}{\lambda}}{\longrightarrow} C\oplus
C\stackrel{\ljz{\pi\e}{\pi}}{\to} M\to 0.$

Note that the homomorphism $\ljz{\pi\e}{\pi}\in\HD((C\oplus
C,\jze),(M,\e))$, so we indeed have an exact sequence of
differential modules.

\vskip 5pt\hskip 70pt$0\to (C\oplus L,\e_{C\oplus
L})\stackrel{\jz{-1}{h}{0}{\lambda}}{\longrightarrow}
\RD{C}\stackrel{\ljz{\pi\e}{\pi}}{\to} (M,\e)\to 0,$\hfill
$(\ddag\ddag)$\vskip 5pt

For any $C'\in\mc{C}$, we have that $\Hr(C',\pi)$ is epi, since
$\Hr(C',\ddag)$ is exact by assumptions. It follows easily that
$\Hr(C',\ljz{\pi\e}{\pi})$ is also epi., which in turn means that
$\Hr(C',\ddag\ddag)$ is exact. Hence we obtain that
$\HD(\RD{C},\ddag\ddag)$ is exact for any $C'\in\mc{C}$, by
Corollary \ref{exC}.

Now take any $R$-module $X$. On one hand, suppose that
$\Hr(\ddag,X)$ is exact. Then, given a homomorphism
$\ljz{g_{_C}}{g_{_L}}\in\Hr(C\oplus L,X)$, there is some
$\theta\in \Hr(C,X)$ such that $g_{_L}=\lambda\theta$. Then we can
check that
$\ljz{g_{_C}}{g_{_L}}=\jz{-1}{h}{0}{\lambda}\ljz{h\theta-g_{_C}}{\theta}$.
Hence $\Hr(\ddag\ddag,X)$ is exact. Consequently,
$\HD(\ddag\ddag,\RD{X})$ is exact by Corollary \ref{exC}.

On the other hand, suppose that $\HD(\ddag\ddag,\RD{X})$ is exact,
which is equivalent to say that  $\Hr(\ddag\ddag,X)$ is exact by
Corollary \ref{exC}.  Given a homomorphism $x\in\Hr(L,X)$, we
obtain $\ljz{0}{x}\in\Hr(C\oplus L,X)$. It follows that there is
some $\ljz{y}{z}\in\Hr(C\oplus C,X)$ such that
$\ljz{0}{x}=\jz{-1}{h}{0}{\lambda}\ljz{y}{z}$. Hence $x=\lambda z$
and consequently $\Hr(\lambda,X)$ is epi. Thus, $\Hr(\ddag,X)$ is
exact.\
\ \hfill$\Box$

%%%%%--------------------------------------------------------------------------

\vskip 15pt

Dually, we have the following result.

\bg{Lem}\label{peL}%preenvelop Lemma

Let $(M,\e)$ be a differential module and $\mc{C}$ be a class of
$R$-modules. Assume that

\vskip 5pt\hskip 90pt$0\to M\stackrel{\lambda}{\to}
C\stackrel{\pi}{\to} L\to 0$\hfill${\ddag^o}$\vskip 5pt

\noindent is an exact sequence of $R$-modules such that
$C\in\mc{C}$ and $\Hr(\ddag^o,C')$ is exact for any $C'\in\mc{C}$.
Then there is an exact sequence of differential modules

\hskip 70pt$0\to
(M,\e)\stackrel{\hjz{\lambda}{\e\lambda}}{\longrightarrow} \RD{C}
\stackrel{\jz{\pi}{h}{0}{-1}}{\longrightarrow} (L\oplus
C,\e_{L\oplus C})\to 0,$\hfill $(\ddag\ddag^o)$

\noindent for some $h\in\Hr(C,C)$, such that

$(1)$ $\HD(\ddag^o,\RD{C'})$ is exact for any $C'\in\mc{C}$.

$(2)$ For any $R$-module $X$, $\Hr(X,\ddag^o)$ is exact if and
only if $\HD(\RD{X},\ddag\ddag^o)$ is exact.

\ed{Lem}

%%%%%--------------------------------------------------------------------------

\vskip 10pt

%The above result shows that $\mc{C}^{\e}$ is a
%precover class in category of differential modules, whenever
%$\mc{C}$ is a precover class in the category of $R$-modules.

Let $\mc{C}$ be a a class of $R$-modules. Recall that a proper
$\mc{C}$-resolution of an $R$-module $M$ is an exact sequence
$\cdots\to C_2\to C_1\to C_0\to M\to 0\ \ (\star)$, where
$C_i\in\mc{C}$ for all $i\ge 0$, such that $\Hr(C,\star)$ is exact
for any $C\in\mc{C}$.

Dually, a proper $\mc{C}$-coresolution of an $R$-module $M$ is an
exact sequence $0\to M\to C_0\to C_1\to C_2\to \cdots\ \
(\star^o)$, where $C_i\in\mc{C}$ for all $i\ge 0$, such that
$\Hr(\star^o,C)$ is exact for any $C\in\mc{C}$.

In the following, we denote by $\mc{C}^{\e}$ the class of
differential modules $\RD{C}$ with
$C\in\mc{C}$.%, where  $\mc{C}$ is a class of $R$-modules. %For instance, Take $\mc{C}$ to be the class of all
%projective $R$-modules, then $\mc{C}^{\e}$ is just the class of
%all projective differential $R$-modules.

\bg{Lem}\label{resL}%resolution lemma

Let $(M,\e_M)\in\D$ and $\mc{C}$ be a class of $R$-modules which
is closed under direct sums. Assume that there is a proper
$\mc{C}$-resolution of the $R$-module $M$, say,

\vskip 5pt\hskip 50pt$\cdots\to C_2\stackrel{c_2}{\to}
C_1\stackrel{c_1}{\to} C_0\stackrel{c_0}{\to} M\to
0$\hfill$(\natural)$\vskip 5pt

\noindent Denote $M_i=\mr{Im}c_i$ for $i\ge 0$. Then there is a
proper $\mc{C}^{\e}$-resolution of the differential module
$(M,\e_M)$

\vskip 5pt$\cdots\to (Q_2\oplus Q_2,\jze)\stackrel{q_2}{\to}
(Q_1\oplus Q_1,\jze)\stackrel{q_1}{\to} (Q_0\oplus
Q_0,\jze)\stackrel{q_0}{\to} (M,\e_M)\to 0$
\hfill$(\natural\natural)$\vskip 5pt

\noindent such that

$(1)$ $Q_i\simeq\oplus_{k=0}^iC_k$ and $\mr{Ker}q_{i}\simeq
Q_i\oplus M_{i+1}$, for all $i\ge 0$;

$(2)$ For any $R$-module $X$, $\Hr(\natural,X)$ is exact if and
only if $\HD(\natural\natural,\RD{X})$ is exact.

\ed{Lem}

\Pf. Note that $M_0=M$ and we have exact sequences

\vskip 5pt\hskip 90pt $0\to M_{i+1}\stackrel{\lambda_i}{\to}
C_i\stackrel{\pi_i}{\to}  M_i\to 0$\hfill$({\natural_i})$\vskip
5pt

\noindent such that $\pi_i\lambda_{i-1}=c_i$ and
$\Hr(C,\natural_i)$ is exact for any $C\in\mc{C}$.

Consider firstly the exact sequence

\vskip 5pt\hskip 90pt $0\to M_1\stackrel{\lambda_0}{\rightarrow}
C_0\stackrel{\pi_0}{\rightarrow} M_0\to
0.$\hfill$({\natural_0})$\vskip 5pt

By Lemma \ref{pcL}, we obtain an exact sequences of differential
modules

\vskip 5pt\hskip 50pt$0\to (C_0\oplus M_1,\e_{C_0\oplus M_1})\to
\RD{C_0}\to (M_0,\e_{M_0})\to
0,$\hfill$(\natural\natural_0)$\vskip 5pt

\noindent such that $\HD(\RD{C},\natural\natural_0)$ is exact  for
any $\RD{C}\in\mc{C}^{\e}$  and such that $\Hr(\natural_0,X)$ is
exact if and only if $\HD(\natural\natural_0,\RD{X})$ is exact,
for any $R$-module $X$.

Note that we have an exact sequence of $R$-modules

\vskip 5pt\hskip 90pt $0\to
M_2\stackrel{\hjz{0}{\lambda_1}}{\longrightarrow} C_0\oplus
C_1\stackrel{\jz{1}{0}{0}{\pi_1}}{\longrightarrow} C_0\oplus
M_1\to 0$\hfill$(\natural_1')$\vskip 5pt

\noindent which comes from $(\natural_1)$ and that
$\Hr(C,\natural_1')$ is exact for any $C\in\mc{C}$, since
$\Hr(C,\natural_1)$ is exact. Thus, by repeating the above process
to the exact sequence $(\natural_1')$ instead of $(\natural_0)$,
and so on, we obtain exact sequences of differential modules

\vskip 5pt\hskip 90pt$0\to (N_{i+1},\e_{N_{i+1}})\to (Q_i\oplus
Q_i,\jze)\to (N_{i},\e_{N_{i}})\to
0,$\hfill$(\natural\natural_i)$\vskip 5pt

\noindent where $Q_i\simeq\oplus_{k=0}^iC_k$ and $N_{i+1}\simeq
Q_i\oplus M_{i+1}$, such that $\HD(\RD{C},\natural\natural_i)$ is
exact for any $\RD{C}\in\mc{C}^{\e}$ and such that
$\Hr(\natural_i',X)$ is exact if and only if
$\HD(\natural\natural_i,\RD{X})$  is exact, for any $R$-module
$X$.

Now the desired sequence $(\natural\natural)$ follows by combining
together the short exact sequences $(\natural\natural_i)$'s.
%%%%%%

\ \ \hfill$\Box$

%%%%%--------------------------------------------------------------------------
\vskip 10pt

We also have the following dual result.

\bg{Lem}\label{coresL}%coresolution lemma

Let $(M,\e_M)\in\D$ and $\mc{C}$ be a class of $R$-modules which
is closed under direct sums. Assume that there is a proper
$\mc{C}$-coresolution of the $R$-module $M$, say,

\vskip 5pt\hskip 50pt$0\to M\stackrel{c_0}{\to}
C_0\stackrel{c_1}{\to} C_1\stackrel{c_2}{\to} C_2\to\cdots\
$\hfill$(\natural^o)$\vskip 5pt

\noindent Denote $M_i=\mr{Im}c_i$ for $i\ge 0$. Then there is a
proper $\mc{C}^{\e}$-coresolution of the differential module
$(M,\e_M)$

\vskip 5pt$0\to(M,\e)\stackrel{q_0}{\to}
\RD{Q_0}\stackrel{q_1}{\to} \RD{Q_1}\stackrel{q_2}{\to}
\RD{Q_2}\to \cdots\ $ \hfill$(\natural\natural^o)$\vskip 5pt

\noindent such that

$(1)$ $Q_i\simeq\oplus_{k=0}^iC_k$ and $\mr{Coker}q_{i}\simeq
M_{i+1}\oplus Q_i$, for all $i\ge 0$;

$(2)$ For any $R$-module $X$, $\Hr(X,\natural^o)$ is exact if and
only if $\HD(\RD{X},\natural\natural^o)$ is exact.

\ed{Lem}

%%%%%--------------------------------------------------------------------------

The following result describes differential modules which are
orthogonal to contractible differential modules.

\bg{Pro}\label{ExtP}%Ext^n Proposition

Let $(M,\e)\in\D$ and $X$ be an $R$-module. Then

$(1)$ $\ED^{i}((M,\e),\RD{X})=0$ for all $i\ge 1$ if and only if
$\Er^i(M,X)=0$ for all $i\ge 0$.

$(2)$ $\ED^{i}(\RD{X},(M,\e))=0$ for all $i\ge 1$ if and only if
$\Er^i(X,M)=0$ for all $i\ge 0$.

\ed{Pro}

\Pf. (1) Take a projective resolution of the $R$-module $M$

\vskip 5pt\hskip 90pt$\cdots\to P_2\stackrel{p_2}{\to}
P_1\stackrel{p_1}{\to} P_0\stackrel{p_0}{\to} M\to
0.$\hfill$(\natural_P)$\vskip 5pt

 By Lemma \ref{resL}, there is a projective resolution of
the differential module $(M,\e_M)$

\vskip 5pt\noindent $\cdots\to (Q_2\oplus
Q_2,\jze)\stackrel{q_2}{\to} (Q_1\oplus
Q_1,\jze)\stackrel{q_1}{\to} (Q_0\oplus
Q_0,\jze)\stackrel{q_0}{\to} (M,\e_M)\to 0,$
\hfill$(\natural\natural_P)$\vskip 5pt

\noindent such that $\Hr(\natural_P,X)$ is exact if and only if
$\HD(\natural\natural_P,\RD{X})$ is exact. Hence, we can  deduce
that

\hskip 90pt $\ED^{i}((M,\e),\RD{X})=0$ for all $i\ge 1$

\hskip 90pt $\Leftrightarrow$ $\HD(\natural\natural_P,\RD{X})$ is
exact

\hskip 90pt $\Leftrightarrow$ $\Hr(\natural_P,X)$ is exact

\hskip 90pt $\Leftrightarrow$ $\Er^i(M,X)=0$ for all $i\ge 1$.

(2). The proof is dual to that of (1).
%%%
\hfill$\Box$

%%%%%--------------------------------------------------------------------------

\vskip 15pt

Given a ring $R$, a Gorenstein projective module $G$ is defined to
be the image of the homomorphism $g_0$ which is given by an exact
sequence of projective $R$-modules

\hskip 90pt $\cdots\to P_1 \stackrel{g_1}{\to} P_0
\stackrel{g_0}{\to} P_{-1} \stackrel{g_{-1}}{\to}\cdots\ $ \hfill
$(\dag^G)$

\noindent such that $\Hr(\dag^G,P)$ is exact for any projective
$R$-module $P$.

Dually,  a Gorenstein injective module $Q$ is defined to be the
image of the homomorphism $q_0$ which is given by an exact
sequence of injective $R$-modules

\hskip 90pt $\cdots\to I_1 \stackrel{q_1}{\to} I_0
\stackrel{q_0}{\to} I_{-1} \stackrel{q_{-1}}{\to}\cdots\ $ \hfill
$(\dag^Q)$

\noindent such that $\Hr(I,\dag^Q)$ is exact for any injective
$R$-module $I$.

%\bg{Sc}\label{PD}\ed{Sc}\vskip -10pt%projective diff mod

Note that the projective differential module over $R$, i.e., the
projective object in $\D$, are just contractible differential
modules such that their underlying modules are projective
$R$-modules [\ref{CE}, Chap. IV]. Similarly, the injective
differential module over $R$, i.e., the injective object in $\D$,
are just contractible differential modules such that their
underlying modules are injective $R$-modules.

We now describe Gorenstein projective (respectively, Gorenstein
injective) differential modules over $R$.

\bg{Th}\label{GT}%Gorenstein Theorem

Let $(M,\e_M)$ be a differential module.

$(1)$ $(M,\e_M)$ is Gorenstein projective if and only if $M$ is a
Gorenstein projective $R$-module.

$(2)$ $(M,\e_M)$ is Gorenstein injective if and only if $M$ is a
Gorenstein injective $R$-module.

\ed{Th}

\Pf. (1) The only-if-part. Since  $(M,\e_M)$ is Gorenstein
projective, there is an exact sequence of projective differential
modules

\hskip 90pt $\cdots\to (P_1,\e_1)\to (P_0,\e_0)\to^f (P_{-1},\e
_{-1})\to\cdots\ $   \hfill$(\ddag^G)$

\noindent such that $(M,\e_M)=\mr{Im}f$ and $\HD(\ddag^G,\RD{P})$
is exact for any projective differential module $\RD{P}$. It
follows from Corollary \ref{exC} that $\Hr(\ddag^G,P)$ is also
exact for any projective module $P$. Hence we see that $M$ is
Gorenstein projective.

The if-part. Assume that $M$ is Gorenstein projective, then there
is an exact sequence of projective modules

\hskip 50pt $\cdots\to P_2\to^{p_2} P_1\to^{p_1}  P_0\to^{p_0}
P_{-1}\to^{p_{-1}} P_{-2}\to^{p_{-2}}\cdots\ $\hfill$(\ast)$

\noindent such that $M=\mr{Im}p_0$ and $\Hr(\ast,P)$ is exact for
any projective module $P$. Then we have two exact sequences

\vskip 5pt\hskip 50pt $\cdots\to P_2\to^{p_2} P_1\to^{p_1}
P_0\to^{\pi_0} M\to 0 $ \hfill$(\ast_l)$\vskip 5pt

\noindent and

\vskip 5pt\hskip 50pt $0\to M\to^{\lambda_0} P_{-1}\to^{p_{-1}}
P_{-2}\to^{p_{-2}}\cdots\ ,$\hfill$(\ast_r)$\vskip 5pt

\noindent such that both $\Hr(\ast_l,P)$ and $\Hr(\ast_r,P)$ are
exact for any projective $R$-module $P$.

By Lemmas \ref{resL} and \ref{coresL}, we obtain two exact
sequences

\vskip 5pt\noindent $\cdots\to \RD{Q_2}\stackrel{q_2}{\to}
\RD{Q_1}\stackrel{q_1}{\to} \RD{Q_0}\stackrel{\pi_0'}{\to}
(M,\e_M)\to 0 $ \hfill$(\ast\ast_l)$\vskip 5pt

\noindent and

\vskip 5pt $0\to (M,\e_M)\stackrel{\lambda_0'}{\to}
\RD{Q_{-1}}\stackrel{q_{-1}}{\to}
\RD{Q_{-2}}\stackrel{q_{-2}}{\to}\cdots\
,$\hfill$(\ast\ast_r)$\vskip 5pt

\noindent where $Q_i$'s are projective $R$-modules, such that both
$\HD(\ast\ast_l,\RD{P})$ and $\Hr(\ast\ast_r,\RD{P}$ are exact for
any projective differential module $\RD{P}$.

Now combining sequences $(\ast\ast_l)$ and $(\ast\ast_r)$ together
we obtain an exact sequence of projective differential modules

\vskip 5pt $\cdots\to \RD{Q_2}\stackrel{q_2}{\to}
\RD{Q_1}\stackrel{q_1}{\to} \RD{Q_0}\stackrel{q_0}{\to} $

\hskip 120pt $\RD{Q_{-1}}\stackrel{q_{-1}}{\to}
\RD{Q_{-2}}\stackrel{q_{-2}}{\to}\cdots\
,$\hfill$(\ast\ast)$\vskip 5pt

\noindent such that $(M,\e_M)=\mr{Im}q_0$ and
$\HD(\ast\ast,\RD{P})$ is exact for any projective differential
module $\RD{P}$. It follows that $(M,\e_M)$ is a Gorenstein
projective differential module.

(2) The proof is dual to that of (1).
%%%%%%%%%%
\hfill$\Box$

%%%%%--------------------------------------------------------------------------

\vskip 15pt

Let $M$ be an $R$-module. The Gorenstein projective dimension of
$M$, denoted by $\mr{Gpd}M$, is defined to be the minimal integer
$n$ such that there is an exact sequence $0\to G_n\to \cdots\to
G_0\to M\to 0$ with all $G_i$'s Gorenstein projective, or $\infty$
if no such exact sequence exists. The Gorenstein injective
dimension of $M$, denoted by $\mr{Gid}M$,  is defined dually.
Moreover, the supresum of Gorenstein projective dimensions of all
$R$-modules coincides with the supresum of Gorenstein injective
dimensions of all $R$-modules, which is called the Gorenstein
global dimension of $R$ and is denoted by $\mr{Ggd}R$ [\ref{BM}].

We have the following result for these Gorenstein homological
dimensions.

\bg{Th}\label{GdT}%Gorenstein dimension Theorem

Let $R$ be a ring and  $(M,\e_M)$ be a differential module and $n$
an integer.

$(1)$ $\mr{Gpd}(M,\e_M)\le n$  if and only if $\mr{Gpd}M\le n$.

$(2)$ $\mr{Gid}(M,\e_M)\le n$  if and only if $\mr{Gid}M\le n$.

$(3)$ $\mr{Ggd}\D\le n$ if and only if $\mr{Ggd}R\le n$.

\ed{Th}

\Pf. (1) Assume first $\mr{Gpd}(M,\e_M)\le n$. Then there is an
exact sequence of differential modules

         \hskip 20pt $0\to (M_n,\e_n)\to \RD{P_{n-1}}\to\cdots\to
         \RD{P_0}\to (M,\e_M)\to 0$

         \noindent such that $P_i$'s are projective $R$-modules
         and $(M_n,\e_n)$ is a Gorenstein projective differential module, see for instance [\ref{Eb}].
Hence, we have an exact sequence of $R$-modules

         \hskip 50pt $0\to M_n\to P_{n-1}\oplus P_{n-1}\to\cdots\to
         P_0\oplus P_0\to M\to 0$.

      \noindent Note that each $P_i\oplus P_i$ is projective and
      that $M_n$ is Gorenstein projective by Theorem \ref{GT}, so
      we have that $\mr{Gpd}M\le n$.

Now assume that $\mr{Gpd}M\le n$. Then we have an exact sequence
of $R$-modules

         \hskip 50pt $0\to L_n\to C_{n-1}\to\cdots\to
         C_0\to M\to 0$

    \noindent such that $L_n$ is Gorenstein projective and that
    each $C_i$ is projective, see also [\ref{Eb}].
By Lemma \ref{resL}, we obtain an exact sequence of differential
modules

         \hskip 20pt $0\to (N_n,\e_n)\to \RD{Q_{n-1}} \to\cdots\to
         \RD{Q_0}\to (M,\e_M)\to 0,$

\noindent such that $Q_i\simeq\oplus_{k=0}^kC_k$ and $N_n\simeq
Q_{n-1}\oplus L_n$. Obviously, $Q_i$'s are projective and $N_n$ is
Gorenstein projective. Hence $(N_n,\e_n)$ is  a Gorenstein
projective differential module. It follows that
$\mr{Gpd}(M,\e_M)\le n$.

(2) Dually to the proof of (1).

(3) By (1) and the definition of Gorenstein global dimension.
%%%%%%%%%%%%%%%%
\hfill$\Box$

%%%%%%%%%%%%%%%%%%%%%%%%%%%

%%%
%%%%%%%%%%%%%%%%%%%%%%%%%%%%%%%%%%%%%%%%%%%%%%%%%%%%%%%%%%%%%%%%%%%%%%%%%%%%%%%%%%%
%%%**************************Ö л**********************************************
%\vskip 20pt
%\begin{center}
%{\bf ACKNOWLEDGEMENTS}
%\end{center}
%
%The author thanks for the referee's carefully reading and helpful
%suggestions.
%%%%%%%%%%%%%%%%%%%%%%%%%%%%%%%%%%%%%%%%%%%%%%%%%%%%%%%%%%%%%%%%%%%%%%%%%%%%%%%%%%
%**************************²Î¿¼ÎÄÏ×**********************************************

\vskip 30pt

{\small

}

\end{document}